\documentstyle[12pt]{article}
\def\thebibliography#1{\bigskip\section*{\centering
References\\}\bigskip\list
{\arabic{enumi}.}{\settowidth\labelwidth{#1}\leftmargin\labelwidth
\advance\leftmargin\labelsep
\usecounter{enumi}}
\def\newblock{\hskip .11em plus .33em minus .07em}
\sloppy\clubpenalty4000\widowpenalty4000
\sfcode`\.=1000\relax}

\def\op#1{\mathop{\fam0 #1}\limits}

\newcommand{\Id}{{\rm Id\,}}

\newcommand{\ben}{\begin{eqnarray}}
\newcommand{\een}{\end{eqnarray}}
\newcommand{\be}{\begin{eqnarray*}}
\newcommand{\ee}{\end{eqnarray*}}
\newcommand{\bea}{\begin{eqalph}}
\newcommand{\eea}{\end{eqalph}}
\newcommand{\cL}{{\cal L}}

\newcommand{\cD}{{\cal D}}
\newcommand{\al}{\alpha}

\newcommand{\la}{\lambda}

\newcommand{\om}{\omega}
\newcommand{\Om}{\Omega}
\newcommand{\m}{\mu}

\newcommand{\g}{\gamma}
\newcommand{\G}{\Gamma}
\newcommand{\e}{\epsilon}

\newcommand{\th}{\theta}

\newcommand{\Si}{\Sigma}
\newcommand{\si}{\sigma}

\newcommand{\w}{\wedge}

\newcommand{\wt}{\widetilde}
\newcommand{\wh}{\widehat}
\newcommand{\ol}{\overline}
\newcommand{\dr}{\partial}

\newcounter{eqalph}
\newcounter{equationa}

\newenvironment{eqalph}{\stepcounter{equation}
\setcounter{equationa}{\value{equation}}
\setcounter{equation}{0}

\begin{eqnarray}}{\end{eqnarray}
\setcounter{equation}{\value{equationa}}}

\newenvironment{proposition}[1]{\medskip{{\sc Proposition}
#1.}}{$\Box$\medskip}
\newenvironment{corollary}[1]{\medskip{{\sc Corollary} #1.}}{$\Box$\medskip}

\newenvironment{remark}{{\bf Remark.}}{$\Box$\medskip}

\begin{document}
\hbox{}

\centerline{\bf\large DIFFERENTIAL GEOMETRY OF COMPOSITE MANIFOLDS.}
\medskip

\centerline{\bf Gennadi A Sardanashvily}
\medskip

\centerline{Department of Theoretical Physics, Physics Faculty,}

\centerline{Moscow State University, 117234 Moscow, Russia}

\centerline{E-mail: sard@grav.phys.msu.su}
\bigskip

\begin{abstract}

In classical field theory, the composite fibred manifolds $Y\to\Si\to X$
provides the adequate mathematical formulation of gauge models with
broken symmetries, e.g., the gauge gravitation theory. This work is
devoted to connections on composite fibred manifolds. In particular, we get
the horizontal splitting of the vertical tangent bundle of a composite
fibred manifold, besides the familiar one of its tangent bundle. This
splitting defines the modified covariant differential and implies the
special fashion of Lagrangian densities of field models on composite manifolds.
The spinor composite bundles are examined.
\end{abstract}

\section{Introduction}

The geometric description of classical
fields by sections of fibred manifolds $Y\to X$ is generally accepted.

\begin{remark}
All maps throughout are of class $C^\infty$. Manifolds are real, Hausdorff,
finite-dimensional, second-countable and connected.
By a fibred manifold is meant a surjective submersion
\[
\pi:Y\to X
\]
provided with an atlas of fibred coordinates $(x^\la, y^i)$.
A locally trivial fibred manifold is termed the bundle.
We denote by $VY$ and $V^*Y$  vertical tangent and
vertical cotangent bundles of $Y$ respectively.
For the sake of simplicity, the pullbacks
\[
Y\op\times_XTX \quad {\rm and} \quad Y\op\times_XT^*X
\]
are denoted by $TX$ and $T^*X$ respectively.
We specify the following types of differential forms on fibred manifolds:

the exterior horizontal forms $ Y\to\op\w^r T^*X$,

the tangent-valued horizontal forms
$Y\to\op\w^r T^*X\op\otimes_Y TY$,
including soldering forms
$Y\to T^*X\op\otimes_YVY$,

and the pullback-valued forms
$Y\to \op\w^r T^*Y\op\otimes_Y TX$ and $Y\to \op\w^r T^*Y\op\otimes_Y T^*X$.
\end{remark}

The present article is devoted to composite fibred manifolds
\begin{equation}
 Y\to \Si\to X \label{I1}
\end{equation}
where $Y\to\Si$ is a bundle denoted by $Y_\Si$ and $\Si\to X$ is a
fibred manifold.

In analytical mechanics, composite fibred manifolds
\[
Y\to\Si\to{\bf R}
\]
characterize systems with variable parameters represented by sections
of $\Si\to{\bf R}$ \cite{6sar}.
In gauge theory of principal connections on a principal bundle $P$
whose structure group is reducible
to its closed subgroup $K$, the composite fibred manifold
\[
P\to P/K\to X
\]
describes spontaneous symmetry breaking
\cite{2sar,sard,9sar}. Global sections of $P/K\to X$ are treated
the Higgs fields.

Application of composite fibred manifolds (\ref{I1}) to field theory is
founded on the following speculations. Given
a global section $h$ of the fibred manifold $\Si\to X$, the restriction $Y_h$
of the bundle $Y_\Sigma$ to $h(X)$ is a fibred submanifold
of $Y\to X$. There is the 1:1 correspondence between
the global sections $s_h$ of $Y_h$ and the global sections of
the composite fibred manifold (\ref{I1}) which cover $h$.
In physical terms, one says that sections $s_h$ of $Y_h$
describe fields in the presence of a background parameter
field $h$, whereas sections
of the composite fibred manifold $Y$ describe all the pairs $(s_h,h)$.
It is important when the bundles $Y_h$ and $Y_{h\neq h'}$ fail to be
equivalent in a sense. The configuration space of these pairs is the
first order jet manifold $J^1Y$ of the composite fibred manifold $Y$ and
their phase space is the Legendre bundle $\Pi$ (\ref{00}) over $Y$.

In particular, the gauge gravitation theory is phrased in terms of
composite fibred manifolds including composite
spinor bundles whose sections describe spinor fields on a world manifold.
\cite{3sar,10sar}.

Let $LX$ be the principal bundle of linear frames in tangent spaces to
a world manifold $X^4$. In gravitation theory, its structure group
\[
GL_4=GL^+(4,{\bf R})
\]
is reduced to the connected Lorentz group $L=SO(3,1)$.
In accordance with the well-known theorem, there is
the 1:1 correspondence between the reduced $L$-principal subbundles $L^hX$ of
$LX$ and the global sections $h$ of the quotient bundle
\begin{equation}
\Si :=LX/L\to X^4. \label{5.15}
\end{equation}
These sections are identified to the tetrad gravitational fields.

Let us consider a bundle of complex Clifford algebras ${\bf C}_{3,1}$
over $X^4$. Its subbundles are both a spinor bundle $S_M\to X^4$ and the
bundle $Y_M\to X^4$ of Minkowski spaces of generating elements of
${\bf C}_{3,1}$. There is the bundle morphism
\[
\g: Y_M\otimes S_M\to S_M
\]
which determines representation of elements of $Y_M$ by Dirac's
$\g$-matrices on elements of the spinor bundle $S_M$. To describe
spinor fields on a world manifold, one requires that the bundle
$Y_M$ is isomorphic to the cotangent bundle $T^*X$ of $X^4$. It takes
place if $Y_M$ is associated with some reduced $L$-principal subbundle $L^hX$
of the linear frame bundle $LX$. Then, there exists the representation
\[
\g_h : T^*X\otimes S_h\to S_h
\]
of cotangent vectors to a world manifold $X^4$ by Dirac's $\g$-matrices
on elements of the spinor bundle $S_h$ associated with the lift of $L^hX$
to a $SL(2,{\bf C})$-principal bundle. Sections of $S_h$ describe
spinor fields in the presence of a tetrad gravitational field $h$.

The crucial point consists in the fact that,
for different sections $h$ and $h'$, the representations $\g_h$ and
$\g_{h'}$ fail to be equivalent. It follows that a
spinor field must be regarded only in a pair with a certain
tetrad field $h$. There is the 1:1 correspondence between
these pairs and the sections of the composite bundle
\begin{equation}
S\to\Si\to X^4 \label{L1}
\end{equation}
where $S\to\Si$ is a spinor bundle associated with the $SL(2,{\bf
C})$-lift of the $L$-principal
bundle $LX\to\Si$.

Dynamics of fields represented by sections of a fibred
manifold $Y\to X$ is phrased in terms of jet manifolds
\cite{car,got,kol,sard,8sar}.

\begin{remark}
Recall that the $k$-order jet manifold $J^kY$ of a fibred
manifold $Y$ comprises the equivalence classes
$j^k_xs$, $x\in X$, of sections $s$ of $Y$ identified by the $(k+1)$
terms of their Taylor series at $x$.
The first order jet manifold $J^1Y$ of
$Y$ is both the fibred manifold $J^1Y\to X$
and the affine bundle $J^1Y\to Y $  modelled on the vector
bundle $T^*X\op\otimes_Y VY$.
It is endowed with the adapted coordinates $(x^\la, y^i, y^i_\la)$:
\[
{y'}^i_\la = (\frac{\dr {y'}^i}{\dr y^j}y_\m^j +
\frac{\dr {y'}^i}{\dr x^\m})\frac{\dr x^\m}{\dr {x'}^\la}.
\]
We identify $J^1Y$ to its image under the  canonical bundle monomorphism
\ben
&&\la:J^1Y\op\to_YT^*X \op\otimes_Y TY,\nonumber \\
&& \la=dx^\la \otimes (\dr_\la + y^i_\la \dr_i). \label{18}
\een
Given a fibred morphism of $\Phi: Y\to Y'$
over a diffeomorphism of $X$, its jet prolongation
$J^1\Phi : J^1Y\to  J^1Y'$ reads
\[
{y'}^i_\mu\circ
J^1\Phi=(\dr_\la\Phi^i+\dr_j\Phi^iy^j_\la)\frac{\dr x^\la}{\dr {x'}^\mu}.
\]
A section $\ol s$ of the fibred jet manifold $J^1Y\to X$ is called
holonomic if it is the jet prolongation $\ol s=J^1s$ of a
section $s$ of $Y\to X$.
There is the 1:1 correspondence between
the global sections
\begin{equation}
\G =dx^\la\otimes(\dr_\la +\G^i_\la\dr_i) \label{M5}
\end{equation}
of the affine jet bundle $J^1Y\to Y$ and the connections on the fibred
manifold $Y\to X$. These global sections form the
affine space modelled on the linear space of soldering forms on $Y$.
Every connection $\G$ on $Y\to X$ yields the first order differential
operator
\be
&& D_\G: J^1Y\op\to_YT^*X\op\otimes_YVY,\\
&&D_\G =(y^i_\la-\G^i_\la)dx^\la\otimes\dr_i,
\ee
on $Y\to X$ which is called the covariant differential relative to the
connection $\G$.
\end{remark}

In firld theory, we can restrict ourselves to the first order
Lagrangian formalism when the configuration space of sections of $Y\to
X$ is the first order jet manifold $J^1Y$ of $Y$. A Lagrangian density on
$J^1Y$ is defined to be a morphism
\be
&& L: J^1Y\to\op\w^nT^*X, \qquad n=\dim X,\\
&&L=\cL\om, \qquad \om=dx^1\w...\w dx^n.
\ee
Note that since the jet bundle $J^1Y\to Y$ is affine, every polynomial
Lagrangian density of field theory factors
\begin{equation}
L: J^1Y\to T^*X\op\otimes_YVY\to \op\w^nT^*X.\label{523}
\end{equation}

The feature of the dynamics of field systems on composite fibred manifolds
consists in the following.

Let $Y$ be a composite manifold (\ref{I1})
provided with the fibred coordinates $(x^\la ,\si^m, y^i)$ where
$( x^\la ,\si^m)$ are fibred coordinates of $\Si\to X$.
Every connection
\begin{equation}
A_\Si=dx^\la\otimes(\dr_\la+\wt A^i_\la\dr_i)
+ d\si^m\otimes(\dr_m+A^i_m\dr_i) \label{S11}
\end{equation}
on the bundle $Y\to\Si$ yields the splitting
\begin{equation}
VY=VY_\Sigma\op\oplus_Y (Y\op\times_\Sigma V\Sigma)\label{146}
\end{equation}
and, as a consequence,
the first order differential operator
\be
&&\wt D: J^1Y\to T^*X\op\otimes_Y VY_\Si, \\
&&\wt D= dx^\la\otimes(y^i_\la-\wt A^i_\la -A^i_m\si^m_\la)\dr_i,
\ee
on $Y$. Let $h$ be a global section
of $\Si\to X$ and $Y_h$ the restriction of the bundle $Y_\Si$ to $h(X)$. The
restriction of $\wt D$ to $J^1Y_h\subset J^1Y$
comes to the familiar covariant differential relative to a certain
connection $A_h$ on $Y_h$.
Thus, it is $\wt D$ that
we may utilize in order to construct a Lagrangian density
\begin{equation}
L: J^1Y\op\to^{\wt D} T^*X\op\otimes_YVY_\Si\to\op\w^nT^*X \label{229}
\end{equation}
for sections of a composite manifold. It should be noted that such a
Lagrangian density is never regular because of the constraint conditions
\[
A^i_m\dr^\mu_i\cL =\dr^\mu_m\cL.
\]
If a Lagrangian density is degenerate, the corresponding
Euler-Lagrange equations are underdetermined.

To describe constraint field systems, the multisymplectic
generalization of the Hamiltonian formalism in mechanics is utilized
\cite{car,gun,6sar,sard,8sar}. In the framework of this approach, the phase
space of sections of $Y\to X$ is the Legendre bundle
\begin{equation}
\Pi=\op\w^n T^*X\op\otimes_Y TX\op\otimes_Y V^*Y \label{00}
\end{equation}
over $Y$. It is provided with the fibred coordinates $(x^\la ,y^i,p^\la_i)$.
Note that
every Lagrangian density  $L$ on $J^1Y$ determines the Legendre morphism
\be
&& \wh L: J^1Y\to \Pi,\\
&& (x^\m,y^i,p^\m_i)\circ\wh L=(x^\m,y^i,\dr^\mu_i\cL).
\ee
Its image plays the role of the Lagrangian constraint space.

The Legendre bundle (\ref{00}) carries the multisymplectic form
\[
\Om =dp^\la_i\w
dy^i\w\om\otimes\dr_\la.
\]
We  say that a connection
$\g$ on the fibred Legendre manifold $\Pi\to X$ is a Hamiltonian
connection if the  form  $\g\rfloor\Om$ is closed.
Then, a muiltimomentum Hamiltonian form $H$ on $\Pi$ is defined to be an
exterior form such that
\[
dH=\g\rfloor\Om
\]
for some Hamiltonian connection $\g$.

The major feature of Hamiltonian systems on composite fibred manifolds
lies in the following \cite{9sar}.
Let $Y$ be a composite fibred manifold (\ref{I1}). The Legendre bundle
$\Pi$ over $Y$ is coordinatized by
\[
(x^\la,\si^m,y^i,p^\la_m,p^\la_i).
\]
Let $A_\Si$ be a connection (\ref{S11}) on the bundle $Y_\Si$.
With a connection $A_\Si$, the splitting
\begin{equation}
\Pi=\op\w^nT^*X\op\otimes_YTX\op\otimes_Y
[V^*Y_\Si\op\oplus_Y (Y\op\times_\Si V^*\Si)]\label{230}
\end{equation}
of the Legendre bundle $\Pi$ is performed as an immediate consequence
of the splitting (\ref{146}).
Given the splitting (\ref{230}), the Legendre
bundle $\Pi$ can be provided with the coordinates
\[
\ol p^\la_i=p^\la_i, \qquad
\ol p^\la_m = p^\la_m +A^i_mp^\la_i
\]
which are compatible with this splitting.

In particular, let $h$ be a global section of the fibred manifold $\Sigma$.
It is readily observed that, given the splitting
(\ref{230}), the submanifold
\[
\{\si=h(x), \ol p^\la_m=0\}
\]
of the Legendre bundle $\Pi$ over $Y$ is isomorphic to
the Legendre bundle $\Pi_h$ over the restriction  $Y_h$ of
$Y_\Si$ to $h(X)$.

Let a multimomentum Hamiltonian form
be associated with a Lagrangian density (\ref{229})
which contains the velocities $\si^m_\mu$ only inside the
differential $\wt D$.
Then, the Lagrangian constraints read
\[
\ol p^\m_m=0. \label{232}
\]

\section{Composite fibred manifolds}

The composite fibred manifold (or simply the composite manifold)
is defined to be composition of surjective submersions
\begin{equation}
\pi_{\Si X}\circ\pi_{Y\Si}:Y\to \Si\to X. \label{1.34}
\end{equation}
Roughly speaking, it is the fibred manifold $Y\to
X$ provided with the particular class of coordinate atlases
$( x^\la ,\si^m,y^i)$:
\[
{x'}^\lambda=f^\lambda(x^\mu), \qquad {\sigma'}^m=f^m(x^\mu,\sigma^n),
\qquad  {y'}^i=f^i(x^\mu,\sigma^n,y^j),
\]
where $(x^\m,\si^m)$ are fibred  coordinates  of
$\Si\to X$ and $y^i$ are bundle
coordinates of $Y_\Si$.
Note that if the fibred manifold $\Si$ also is a bundle, the
composite manifold $Y$ fails to be a bundle in general. We further
propose that $\Si$ has a global section.

Recall the following assertions \cite{sard,sau}.

\begin{proposition}{1} Let $Y$ be the composite manifold
(\ref{1.34}). Given a section $h$ of
$\Sigma$ and a section $s_\Sigma$ of $Y_\Sigma$, their
composition  $s_\Sigma\circ h$ obviously
is a section of the composite manifold $Y$. Conversely, if the bundle
$Y_\Si$ has a global section, every
global section $s$ of the fibred manifold $Y\to X$ is  represented by some
composition $s_\Si\circ h$ where $h=\pi_{Y\Si}\circ s$ and $s_\Si$ is an
extension of the local section $h(X)\to s(X)$ of the bundle
$Y_\Si$ over the closed imbedded submanifold $h(X)\subset\Si$.
\end{proposition}

\begin{proposition}{2}
Given a global section $h$ of $\Sigma$, the
restriction $Y_h=h^*Y_\Si$
of the bundle $Y_\Sigma$ to $h(X)$ is a fibred imbedded submanifold
of $Y$.
\end{proposition}

\begin{corollary}{3}
There is the 1:1 correspondence between the
sections $s_h$ of $Y_h$ and the sections $s$ of
the composite manifold $Y$ which cover $h$.
\end{corollary}

Given fibred
coordinates $(x^\la, \si^m, y^i)$ of the composite manifold $Y$,
the jet manifolds $J^1\Si$,
$J^1Y_\Si$ and $J^1Y$ are coordinatized respectively by
\[
( x^\la ,\si^m, \si^m_\la),\qquad
( x^\la ,\si^m, y^i, \wt y^i_\la, y^i_m),\qquad
( x^\la ,\si^m, y^i, \si^m_\la ,y^i_\la).
\]

\begin{proposition}{4} There exists  the  canonical surjection
\ben
&&\rho : J^1\Si\op\times_\Si J^1Y_\Sigma\to J^1Y,\nonumber \\
&&\rho(j^1_xh,j^1_{h(x)}s_\Sigma)=j^1_x({s_\Sigma}\circ h), \label{1.38}\\
&&y^i_\lambda\circ\rho=y^i_m{\sigma}^m_{\lambda} +\wt y^i_{\lambda},\nonumber
\een
where $s_\Si$ and $h$ are sections of $Y_\Si$ and $\Si$ respectively.
\end{proposition}

\begin{corollary}{5}
Let $A_\Si$ be the connection (\ref{S11}) on the bundle $Y_\Si$ and
$\G$ the connection (\ref{M5}) on the fibred manifold $\Sigma$.
Building on the morphism (\ref{1.38}), one can construct the composite
connection
\begin{equation}
A=dx^\la\otimes[\dr_\la+\G^m_\la\dr_m +(A^i_m\G^m_\la + \wt A^i_\la)
\dr_i] \label{1.39}
\end{equation}
on the composite manifold $Y$ in accordance with the commutative
diagram
\[
\begin{array}{rcccl}
 & {J^1\Si\op\times_\Si J^1Y_\Si} &  \op\longrightarrow^{\rho} & {J^1Y} &  \\
{_{\G\times A_\Si}} &\put(0,-10){\vector(0,1){20}} & &
\put(0,-10){\vector(0,1){20}} & {_A} \\
 & {\Si\op\times_\Si Y_\Si} & \op\longleftarrow_{\pi_{Y\Si}\times\Id_Y}
& {Y} &
\end{array}
\]
\end{corollary}

Let a global section $h$ of
$\Si$ be an integral section of the connection $\Gamma$ on
$\Si$, that is, $\Gamma\circ h=J^1h$.
Then, the composite connection (\ref{1.39}) on $Y$ is reducible to the
connection
\begin{equation}
A_h=dx^\lambda\otimes[\dr_\lambda+(A^i_m\dr_\lambda h^m
+\wt A^i_\lambda)\dr_i] \label{1.42}
\end{equation}
on the fibred submanifold $Y_h$ of $Y\to X$  in accordance
with the commutative diagram
\[
\begin{array}{rcccl}
 & {J^1Y_h} &  \op\longrightarrow^{J^1i_h} & {J^1Y} &  \\
{_{A_h}} &\put(0,-10){\vector(0,1){20}} & & \put(0,-10){\vector(0,1){20}} &
{_A} \\
 & {Y_h} & \op\hookrightarrow_{i_h} & {Y} &
\end{array}
\]
In particular, every connection $A_\Si$ (\ref{S11}) on $Y_\Si$,
whenever $h$, is reducible to the connection (\ref{1.42}) on $Y_h$.

\section{Connections on composite manifolds}

Given a composite  manifold (\ref{1.34}), we have the exact sequences
\be
&& 0\to VY_\Si\hookrightarrow VY\to Y\op\times_\Si V\Si\to 0,\\
&& 0\to Y\op\times_\Si V^*\Si\hookrightarrow V^*Y\to V^*Y_\Si\to 0
\ee
over $Y$, besides the familiar ones
\be
&& 0\to VY\hookrightarrow TY\to Y\op\times_X TX\to 0, \\
&& 0\to Y\op\times_X T^*X\hookrightarrow T^*Y\to V^*Y\to 0.
\ee

\begin{proposition}{6}
There exist the canonical splittings
\[
J^1Y_\Si\op\times_YVY=VY_\Si\op\oplus_{J^1Y} (Y\op\times_\Si V\Si),
\]
\[
\dot y^i\dr_i + \dot\si^m\dr_m=
(\dot y^i -y^i_m\dot\si^m)\dr_i + \dot\si^m(\dr_m+y^i_m\dr_i),
\]
of the vertical tangent bundle $VY$ of $Y\to X$ and
\[
J^1Y\op\times_YV^*Y=V^*Y_\Si\op\oplus_{J^1Y} (Y\op\times_\Si V^*\Si),
\]
\[
\dot y_i dy^i + \dot\si_m d\si^m=
\dot y_i(dy^i -y^i_m d\si^m) + (\dot\si_m +y^i_m\dot y_i) d\si^m,
\]
of the vertical cotangent bundle $V^*Y$ of $Y$.
\end{proposition}

The proof is starightforward.

These splittings add the familiar canonical horizontal splittings
\be
&&J^1Y\op\times_Y TY=TX\op\oplus_{J^1Y} VY,\\
&&J^1Y\op\times_Y T^*Y=T^*X\op\oplus_{J^1Y}V^*Y.
\ee

\begin{corollary}{7}
Every connection (\ref{S11}) on the bundle $Y_\Si$ determines:
\begin{itemize}\begin{enumerate}
\item the horizontal splitting
\begin{equation}
VY=VY_\Si\op\oplus_Y (Y\op\times_\Si V\Si),\label{46}
\end{equation}
\[
\dot y^i\dr_i + \dot\si^m\dr_m=
(\dot y^i -A^i_m\dot\si^m)\dr_i + \dot\si^m(\dr_m+A^i_m\dr_i),
\]
of the vertical tangent bundle $VY$ of $Y$,
\item the dual horizontal splitting
\begin{equation}
V^*Y=V^*Y_\Si\op\oplus_Y (Y\op\times_\Si V^*\Si),\label{46'}
\end{equation}
\[
\dot y_i dy^i + \dot\si_m d\si^m=
\dot y_i(dy^i -A^i_m d\si^m) + (\dot\si_m +A^i_m\dot y_i) d\si^m,
\]
of the vertical cotangent bundle $V^*Y$ of $Y$.
\end{enumerate}\end{itemize}
\end{corollary}

It is readily observed that the splittings (\ref{46}) and (\ref{46'})
are uniquely characterized by the form
\begin{equation}
\om\w A_\Si = \om\w d\si^m\otimes(\dr_m +A^i_m\dr_i),\label{227}
\end{equation}
and different connections $A_\Si$ can define the same splittings (\ref{46})
and (\ref{46'}).

Building on the horizontal splitting (\ref{46}), one can constract
the following first order differential operator on the composite manifold
$Y$:
\ben
&&\wt D={\rm pr}_1\circ D_A
: J^1Y\to T^*X\op\otimes_Y VY \to T^*X\op\otimes_Y VY_\Si,
\nonumber\\
 &&\wt D= dx^\la\otimes[y^i_\la-A^i_\la
-A^i_m(\si^m_\la-\G^m_\la)]\dr_i =
 dx^\la\otimes(y^i_\la-\wt
A^i_\la -A^i_m\si^m_\la)\dr_i, \label{7.10}
 \een
where $D_A$ is the covariant differential relative to the
composite connection $A$ (\ref{1.39})
which is composition of $A_\Si$ and some
connection $\G$ on $\Si$. We shall call $\wt D$ the vertical covariant
differential. This possesses the following property.

Given a global section $h$ of $\Si$, let $\G$ be a connection on $\Si$
whose integral section is $h$.
It is readily observed that the vertical covariant differential
(\ref{7.10}) restricted to $J^1Y_h\subset J^1Y$
comes to the familiar covariant
differential relative to the connection $A_h$ (\ref{1.42}) on $Y_h$.
Thus, it is the vertical covariant differential (\ref{7.10}) that
we may utilize in order to construct a Lagrangian density (\ref{229})
for sections of a composite manifold.

Now, we consider connections on a composite manifold
$Y$ when $Y_\Si$
is a vector bundle. Let a connection $A$ on $Y$ be  projectable to a
connection $\G$ on $\Si$ in accordance with the commutative diagram
\[
\begin{array}{rcccl}
 & {J^1Y} &  \op\longrightarrow^{J^1\pi_{Y\Si}} & {J^1\Si} &  \\
{_A} &\put(0,-10){\vector(0,1){20}} & & \put(0,-10){\vector(0,1){20}} & {_\G}
\\
 & {Y} & \op\longrightarrow_{\pi_{Y\Si}} & {\Si} &
\end{array}
\]

Let
\begin{equation}
A=dx^\lambda\otimes(\dr_\lambda+\Gamma^m_\lambda (\si)
\dr_m + A^i{}_{j\lambda}(\si)y^j\dr_i) \label{520}
\end{equation}
be a linear morphism over $\G$. The following
constructions generalize the familiar notions of a dual linear
connection and a tensor product linear connection on vector bundles.

 Let $Y^*\to\Si\to X$ be a composite manifold where $Y^*\to \Si$
is the vector bundle dual to $Y_\Si$. Given the projectable connection
(\ref{520}) on $Y$ over $\G$, there exists the unique projectable connection
\[
A^*=dx^\lambda\otimes(\dr_\lambda+\Gamma^m_\lambda
\dr_m - A^j{}_{i\lambda}y_j\dr^i).
\]
on $Y^*\to X$ over $\G$ such that the following diagram commutes:
\[
\begin{array}{rcccl}
 & {J^1Y\op\times_{J^1\Si} J^1Y^*} & \op\longrightarrow^{J^1\langle\rangle} &
{J^1\Si\op\times_\Si (T^*X\times{\bf R})} & \\
{_{A\times A^*}} &\put(0,-10){\vector(0,1){20}} & &
\put(0,-10){\vector(0,1){20}} & {_{\G\times\wh 0\times\Id}} \\
 & {Y\op\times_\Si Y^*} &  \op\longrightarrow_{\langle\rangle} &
{\Si\times{\bf R}} &
 \end{array}
\]
where $\wh 0$ is the zero section of $T^*X$.
We term $A^*$ the connection dual to $A$ over $\G$.

 Let $Y\to\Si\to X$ and $Y'\to\Si\to X$ be composite manifolds
where $Y\to\Si$ and $Y'\to\Si$ are vector bundles. Let
$A$ and $A'$ be the connections (\ref{520}) on $Y$ and $Y'$
respectively which are projectable to the same connection $\G$ on
$\Si$. There is the unique projectable  connection
\begin{equation}
A\otimes A'=dx^\lambda\otimes[\dr_\lambda+\Gamma^m_\lambda
\dr_m+(A^i{}_{j\la}y^{jk}+{A'}^k{}_{j\la}y^{ij})\dr_{ik}]\label{66}
\end{equation}
 on the tensor product $Y\op\otimes_\Si Y'\to X$
such that the following diagram commutes:
\[
\begin{array}{rcccl}
 & {J^1Y\op\times_{J^1\Si} J^1Y'} & \op\longrightarrow^{J^1\otimes} &
{J^1(Y\op\otimes_{J^1\Si} Y')} &  \\
{_{A\times A'}} &\put(0,-10){\vector(0,1){20}} & &
\put(0,-10){\vector(0,1){20}} & {_{A\otimes A'}} \\
 & {Y\op\times_{\Si} Y'} &  \op\longrightarrow_{\otimes} &
{Y\op\otimes_{\Si} Y'} &
\end{array}
\]
It is called the tensor product connection over $\G$.

In particular, let $Y\to X$ be a fibred manifold and
 $\G$  a connection on  $Y$. The vertical tangent
morphism $V\G$ to $\G$ determines the connection
\be
&&  V\G :VY\to VJ^1Y=J^1VY, \\
&& V\G =dx^\la\otimes(\dr_\la
+\G^i_\la\frac{\dr}{\dr y^i}+\dr_j\G^i_\la\dot y^j \frac{\dr}{\dr \dot y^i}),
 \ee
on the composite manifold $VY\to Y\to X$.
 The connection $V\G$ is projectable to $\G$,
and it is a linear bundle morphism over $\G$.
It yields the  connection
\begin{equation}
V^*\G =dx^\la\otimes(\dr_\la +\G^i_\la\frac{\dr}{\dr
y^i}-\dr_j\G^i_\la \dot y_i \frac{\dr}{\dr \dot y_j}) \label{44}
\end{equation}
on the composite manifold $ V^*Y\to Y\to X$ which is the dual connection to
$V\G$ over $\G$.

For instance, every connection $\G$ on a fibred manifold $Y\to X$
gives rise to the connection
\begin{equation}
\wt\G =dx^\la\otimes[\dr_\la +\G^i_\la (y)\dr_i +
(-\dr_j\G^i_\la (y)  p^\m_i-K^\m{}_{\nu\la}(x) p^\nu_j+K^\al{}_{\al\la}(x)
p^\m_j)\dr^j_\m]  \label{404}
\end{equation}
on the fibred Legendre manifold $\Pi\to X$ where $K$
is a symmetric linear connection
\be
&& K^\al_\la=-K^\al{}_{\nu\la}(x)\dot x^\nu,\\
&&K^*_{\al\la}=K^\nu{}_{\al\la}(x)\dot x_\nu, \\
&& K^\m{}_{\nu\la}=K^\m{}_{\la\nu},
\ee
on the bundles $TX$ and $T^*X$.
The connection (\ref{404}) is the tensor product (\ref{66}) [over $\G$] of
the connection $\G\times K$ on the pullback composite manifold
\[
Y\op\times_X\op\w^{n-1}T^*X\to Y\to X
\]
and the connection $V^*Y$ (\ref{44}) on the composite manifold $V^*Y\to
Y\to X$. The connection (\ref{404}) covers the
connection $\G$ on the fibred manifold $Y\to X$.
These connections play the prominent role in the multisymplectic
Hamiltonian formalism.
Let $\G$ be a connection on the fibred manifold $Y\to X$ and
\be
&&\wh\G:=\G\circ\pi_{\Pi Y}:\Pi\to Y\to J^1Y, \\
&&\wh\G=dx^\la\otimes (\dr_\la +\G^i_\la\dr_i),
\ee
its pullback by $\pi_{\Pi Y}$ onto the Legendre bundle $\Pi$ over $Y$.
Then, the lift $\wt\G$ (\ref{404}) of $\G$ onto
the fibred Legendre manifold $\Pi\to X$ obeys the identity
\[
\wt\G\rfloor\Om =d(\wh\G\rfloor\th).
\]
A glance at this identity shows that $\wt\G$ is a Hamiltonian connection.

\section{Composite fibration of principal bundles}

Let $\pi_P:P\to X$ be a principal bundle with a structure
Lie group $G$ which acts freely and
transitively on $P$ on the right:
\[
r_g : p\mapsto pg, \qquad
p\in P,\quad g\in G.
\]

Let $K$ be a closed subgroup of $G$. We have the composite manifold
\begin{equation}
\pi_{\Si X}\circ\pi_{P\Si}: P\to P/K\to X \label{7.16}
\end{equation}
where
\[
P_\Si:=P\to P/K
\]
is a principal bundle with the structure group $K$ and
\[
\Si_K=P/K=(P\times G/K)/G
\]
is the $P$-associated bundle with the standard fiber $G/K$ on which the
structure group $G$ acts on the left.

Let the structure group $G$ be reducible to its closed subgroup $K$.
Recall the 1:1 correspondence
\[
\pi_{P\Si}(P_h)=(h\circ\pi_P)(P_h)
\]
between the global sections $h$ of the
bundle $P/K\to X$ and the reduced $K$-principal
subbundles $P_h$ of $P$ which consist with
restrictions of the principal bundle $P_\Si$ to $h(X)$.

Given the composite manifold (\ref{7.16}), the
canonical morphism (\ref{1.38}) results in the surjection
\[J^1P_\Si/K\op\times_\Si J^1\Si \to J^1P/K \]
over $J^1\Si$. Let $A_\Si$ be a principal connection on $P_\Si$
and $\G$ a connection on $\Si$. The corresponding composite connection
(\ref{1.39}) on the composite manifold (\ref{7.16})
is equivariant under the canonical
action of $K$ on $P$. If the connection $\G$ has an integral global section
$h$ of $P/K\to X$, the composite connection (\ref{1.39}) is reducible to the
connection (\ref{1.42}) on $P_h$ which consists with the principal connection
on $P_h$ induced by $A_\Si$.

Let us consider the composite manifold
\begin{equation}
Y=(P\times V)/K\to P/K\to X \label{7.19}
\end{equation}
where the bundle
\[
Y_\Si:= (P\times V)/K\to P/K
\]
is associated with the $K$-principal bundle
$P_\Si$. Given a reduced subbundle $P_h$ of $P$,
the associated bundle
\[
Y_h=(P_h\times V)/K
\]
is isomorphic to the restriction of $Y_\Si$ to $h(X)\subset \Si_K$.

Note that the manifold $(P\times V)/K$ possesses also the structure of the
bundle
\[
Y=(P\times (G\times V)/K)/G
\]
associated with the principal bundle $P$. Its standard fibre
is $(G\times V)/K$ on which the structure group $G$ of $P$ (and its
subgroup $K$) acts by the law
\[
G\ni g: (G\times V)/K\to (gG\times V)/K.
\]
It differs from action of the structure group $K$ of $P_\Si$ on this
standard fibre. As a
shorthand, we can write the latter in the form
\[
K\ni g: (G\times V)/K\to (G\times gV)/K.
\]
However, this action fails to be canonical and depends on existence and
specification of a global section of the bundle $G\to G/K$.
If the standard fibre $V$ of the bundle $Y_\Si$ carriers representation
of the whole group $G$, these two actions are equivalent, otherwise
in general case.

\section{Composite spinor bundles}

By $X^4$ is further meant an oriented world manifold which satisfies
the well-known global topological conditions in order that a spinor
structure can exist. To
summarize these conditions, we assume that $X^4$ is not compact and
the linear frame bundle $LX$ over $X^4$ is trivial.

Given
a Minkowski space $M$ with the Minkowski metric $\eta$, let
\[
A_M=\op\oplus_nM^n,\qquad M^0={\bf R},\qquad
M^{n>0}=\op\otimes^nM,
\]
be the tensor algebra modelled on $M$. The complexified quotient of this
algebra by the two-sided ideal generated by elements
\[
e\otimes e'+e'\otimes e-2\eta(e,e')\in A_M,\qquad e\in M,
\]
constitutes the complex Clifford algebra ${\bf C}_{1,3}$.
A spinor space $V$ is defined to be a
linear space of some minimal left ideal of ${\bf C}_{1,3}$  on
which this algebra acts on the left. We have the representation
\begin{equation}
\g: M\otimes V \to V \label{521}
\end{equation}
of elements of the Minkowski space $M\subset{\bf C}_{1,3}$ by
Dirac's $\g$-matrices on $V$:
\[
\g (e^a\otimes y^Av_A) = \g^{aA}{}_By^Bv_A,
\]
where
$\{e^0...e^3\}$ is a fixed basis for $M$, $v_A$ is a basis for $V$.

Let us consider the transformations preserving the representation
(\ref{521}).
These are pairs $(l,l_s)$ of Lorentz transformations $l$ of the Minkowski
space $M$ and invertible elements $l_s$ of ${\bf C}_{1,3}$ such that
\be
&&lM = l_s Ml^{-1}_s,\\
&&\g (lM\otimes l_sV) = l_s\g (M\otimes V).
\ee
Elements $l_s$  form the Clifford group whose action on $M$
however is not effective. We restrict ourselves to its spinor
subgroup $L_s =SL(2,{\bf C})$
whose generators act on $V$ by the representation
\[
I_{ab}=\frac{1}{4}[\g_a,\g_b].
\]

Let us consider a bundle of complex Clifford algebras ${\bf C}_{3,1}$
over $X^4$. Its subbundles are both a spinor bundle $S_M\to X^4$ and the
bundle $Y_M\to X^4$ of Minkowski spaces of generating elements of
${\bf C}_{3,1}$.
To describe spinor fields on a world manifold, one must
require $Y_M$ be isomorphic to the cotangent bundle $T^*X$
of a world manifold $X^4$. It
takes place if the linear frame bundle $LX$
contains a reduced $L$ subbundle $L^hX$ such that
\[
Y_M=(L^hX\times M)/L.
\]
In this case, the spinor bundle $S_M$ is associated with the $L_s$-lift
$P_h$ of $L^hX$:
\begin{equation}
S_M=S_h=(P_h\times V)/L_s.\label{510}
\end{equation}

There is the above-mentioned 1:1 correspondence between the reduced
subbubdles $L^hX$ of $LX$ and
global sections $h$ of the bundle $\Si$ (\ref{5.15}).

Given a tetrad field $h$, let $\Psi^h$ be an atlas of
$LX$ such that the corresponding local
sections $z_\xi^h$ of $LX$ take their values into $L^hX$.
With respect to $\Psi^h$ and a
holonomic atlas $\Psi^T=\{\psi_\xi^T\}$ of $LX$, a tetrad field $h$
can be represented by a family of $GL_4$-valued tetrad functions
\[
h_\xi=\psi^T_\xi\circ z^h_\xi.
\]
They are $GL_4$-valued functions of atlas transformations
\begin{equation}
dx^\la= h^\la_a(x)h^a \label{L6}
\end{equation}
between the fibre bases $\{dx^\la\}$ and $\{h^a\}$ for
$T^*X$ associated with $\Psi^T$ and $\Psi^h$ respectively.

Given a tetrad field $h$, one can define the representation
\begin{equation}
\gamma_h: T^*X\otimes S_h= (P_h\times (M\otimes V))/L_s\to (P_h\times
\gamma(M\otimes V))/L_s = S_h \label{L4}
\end{equation}
of cotangent vectors to a world manifold $X^4$ by Dirac's $\g$-matrices
on elements of the spinor bundle $S_h$. With respect to an atlas
$\{z_\xi\}$ of $P_h$ and the associated atlas $\{z^h_\xi\}$ of $LX$,
the morphism (\ref{L4}) reads
\[
\g_h(h^a\otimes y^Av_A(x))=\g^{aA}{}_By^Bv_A(x)
\]
where $\{v_A(x)\}$ are the associated fibre bases
for $S_h$. As a shorthand, one can write
\[
\wh dx^\la=\g_h(dx^\la)=h^\la_a(x)\g^a.
\]

On the physical level, one can say
that, given the representation (\ref{L4}), sections of
the spinor bundle $S_h$ describe spinor fields in the presence of
the tetrad gravitational field $h$. Indeed,
let $A_h$ be a principal connection on $S_h$ and
\be
&&D: J^1S_h\op\to_{S_h}T^*X\op\otimes_{S_h}VS_h,\\
&&D=(y^A_\la-A^{ab}{}_\la (x)I_{ab}{}^A{}_By^B)dx^\la\otimes\dr_A,
\ee
the corresponding covariant differential. Given the
representation (\ref{L4}), one can construct the Dirac operator
\begin{equation}
 \cD_h=\g_h\circ D: J^1S_h\to T^*X\op\otimes_{S_h}VS_h\to VS_h, \label{I13}
\end{equation}
\[
\dot y^A\circ\cD_h=h^\la_a(x)\g^{aA}{}_B(y^B_\la-A^{ab}{}_\la
I_{ab}{}^A{}_By^B).
\]
We here use the fact that the vertical tangent bundle $VS_h$ admits the
canonical splitting
\[
VS_h=S_h\times S_h,
\]
and $\g_h$ in the expression (\ref{I13}) is the pullback
\be
&& \g_h: T^*X\op\otimes_{S_h}VS_h\op\to_{S_h}VS_h,\\
&& \g_h(h^a\otimes\dot y^A\dr_A)=\g^{aA}{}_B\dot y^B\dr_B,
\ee
over $S_h$ of the bundle morphism (\ref{L4}).

In the presence of different tetrad fields $h$ and $h'$, spinor fields are
described by sections of spinor bundles $S_h$ and $S_{h'}$ associated
with $L_s$-lifts $P_h$ and $P_{h'}$ of different reduced $L$-principal
subbundles
of $LX$. Therefore, the representations $\g_h$ and $\g_{h'}$
(\ref{L4}) are not equivalent \cite{3sar,10sar}.
It follows that a spinor field must be regarded only in a pair with
a certain tetrad field. In accordance with Corollary 3,
there is the 1:1 correspondence
between these pairs and sections of the composite spinor bundle (\ref{L1}).

We have the composite manifold
\begin{equation}
\pi_{\Si X}\circ\pi_{P\Si}:LX\to\Si\to X^4 \label{L3}
\end{equation}
where $\Si$ is the quotient bundle (\ref{5.15}) and
\[
LX_\Si:=LX\to\Si
\]
is the L-principal bundle.

Building on the double universal covering of the group $GL_4$, one can
perform the $L_s$-principal lift $P_\Si$ of $LX_\Si$ such that
\[
P_\Si/L_s=\Si, \qquad LX_\Si=r(P_\Si).
\]
In particular, there is imbedding of $P_h$ onto the restriction of $P_\Si$
to $h(X^4)$.

Let us consider the composite spinor bundle (\ref{L1}) where
\[
S_\Si= (P_\Si\times V)/L_s
\]
is
associated with the $L_s$-principal bundle $P_\Si$. It is readily observed
that, given a global section $h$ of $\Si$, the restriction $S_\Si$ to
$h(X^4)$ is the spinor bundle $S_h$ (\ref{510}) whose sections describe
spinor fields in the presence of the tetrad field $h$.

Let us provide the principal bundle $LX$ with a holonomic atlas
$\{\psi^T_\xi, U_\xi\}$ and the principal bundles $P_\Si$ and $LX_\Si$
with associated atlases $\{z^s_\e, U_\e\}$ and $\{z_\e=r\circ z^s_\e\}$.
With respect to these atlases, the composite spinor bundle is endowed
with the fibred coordinates $(x^\la,\si_a^\mu, y^A)$ where $(x^\la,
\si_a^\mu)$ are fibred coordinates of the bundle $\Si$ such that
$\si^\mu_a$ are the matrix components of the group element
\[
GL_4\ni (\psi^T_\xi\circ z_\e)(\si): {\bf R}^4\to {\bf R}^4,
\qquad \si\in U_\e,\qquad \pi_{\Si X}(\si)\in U_\xi.
\]
Given a section $h$ of $\Si$, we have
\be
&&z^h_\xi (x)= (z_\e\circ h)(x), \qquad h(x)\in U_\e,
\qquad x\in U_\xi,\\
&& (\si^\la_a\circ h)(x)= h^\la_a(x),
\ee
where $h^\la_a(x)$ are tetrad functions (\ref{L6}).

The jet manifolds $J^1\Si$, $J^1S_\Si$ and $J^1S$ are coordinatized by
\[
(x^\la,\si^\mu_a, \si^\mu_{a\la}),\qquad
(x^\la,\si^\mu_a, y^A,\wt y^A_\la, y^A{}^a_\mu),\qquad
(x^\la,\si^\mu_a, y^A,\si^\mu_{a\la}, y^A_\la).
\]
Note that, whenever $h$, the jet manifold
$J^1S_h$ is a fibred submanifold of $J^1S\to X^4$ given by the coordinate
relations
\[
\si^\mu_a=h^\mu_a(x), \qquad \si^\mu_{a\la}=\dr_\la h^\mu_a(x).
\]

Let us consider the bundle of Minkowski spaces
\[
(LX\times M)/L\to\Si
\]
associated with the $L$-principal bundle $LX_\Si$. Since $LX_\Si$ is
trivial, it is isomorphic to
the pullback $\Si\op\times_X T^*X$ which we denote by the same symbol
$T^*X$. Building on the morphism (\ref{521}), one can define the bundle
morphism
\begin{equation}
\g_\Si: T^*X\op\otimes_\Si S_\Si= (P_\Si\times (M\otimes V))/L_s
\to (P_\Si\times\g(M\otimes V))/L_s=S_\Si, \label{L7}
\end{equation}
\[
\wh dx^\la=\g_\Si (dx^\la) =\si^\la_a\g^a,
\]
over $\Si$. When restricted to $h(X^4)\subset \Si$,
the morphism (\ref{L7}) comes to the morphism $\g_h$
(\ref{L4}). Because of the canonical vertical splitting
\[
VS_\Si =S_\Si\op\times_\Si S_\Si,
\]
the morphism (\ref{L7}) yields the corresponding morphism
\begin{equation}
\g_\Si: T^*X\op\otimes_SVS_\Si\to VS_\Si. \label{L8}
\end{equation}

We use this morphism in order to construct the total Dirac
operator on sections of the composite spinor bundle $S$ (\ref{L1}). We are
based on the following fact.

Let
\begin{equation}
\wt A=dx^\la\otimes (\dr_\la +\wt A^B_\la\dr_B) + d\si^\mu_a\otimes
(\dr^a_\mu+A^B{}^a_\mu\dr_B) \label{200}
\end{equation}
be a connection on the bundle $S_\Si$. It determines the horizontal
splitting (\ref{46}) of the vertical tangent bundle $VS$ and the
vertical covariant differential (\ref{7.10}).
The composition of
morphisms (\ref{L8}) and (\ref{7.10}) is the first order differential
operator
\[
\cD=\g_\Si\circ\wt D:J^1S\to T^*X\op\otimes_SVS_\Si\to VS_\Si,
\]
\[
\dot y^A\circ\cD=\si^\la_a\g^{aA}{}_B(y^B_\la-\wt A^B_\la -
A^B{}^a_\mu\si^\mu_{a\la}),
\]
on $S$.
One can treat it the total Dirac operator since, whenever a
tetrad field $h$, the restriction of $\cD$ to $J^1S_h\subset J^1S$ comes
to the Dirac operator $\cD_h$ (\ref{I13})
with respect to the connection
\[
A_h=dx^\la\otimes[\dr_\la+(\wt A^B_\la+A^B{}^a_\mu\dr_\la h^\mu_a)\dr_B]
\]
on $S_h$.

To construct the connection (\ref{200}) in explicit form, let us set up the
principal connection onthe bundle $LX_\Si$ which is given by
the local connection form
\ben
&& A_\Si = (\wt A^{ab}{}_\mu dx^\mu+ A^{ab}{}^c_\mu d\si^\mu_c)\otimes I_{ab},
\label{L10}\\
&&\wt A^{ab}{}_\mu=\frac12 K^\nu{}_{\la\mu}\si^\la_c (\eta^{cb}\si^a_\nu
-\eta^{ca}\si^b_\nu ),\nonumber\\
&&A^{ab}{}^c_\mu=\frac12(\eta^{cb}\si^a_\mu -\eta^{ca}\si^b_\mu),
\label{M4}
\een
where $K$ is some symmetric connection on $TX$ and (\ref{M4})
corresponds to the canonical left-invariant free-curvature connection on
the bundle
\[
GL_4\to GL_4/L.
\]
Given a tetrad field $h$, the connection (\ref{L10}) is reduced to the
Levi-Civita connection
\[
A_h =\frac12[K^\nu{}_{\la\mu}\si^\la_c (\eta^{cb}\si^a_\nu
-\eta^{ca}\si^b_\nu)+\dr_\mu h^\nu_c(\eta^{cb}\si^a_\nu-\eta^{ca}\si^b_\nu)]
\]
on $L^hX$.

The connection (\ref{200}) on the spinor bundle $S_\Si$ which is
associated with $A_\Si$ (\ref{L10}) reads
\[
\wt A_\Si=dx^\la\otimes (\dr_\la +\frac12\wt A^{ab}{}_\la
I_{ab}{}^B{}_Ay^A\dr_B) + d\si^\mu_c\otimes
(\dr^c_\mu+\frac12 A^{ab}{}^c_\mu I_{ab}{}^B{}_Ay^A\dr_B).
\]
It determines the canonical horizontal
splitting (\ref{46}) of the vertical tangent bundle $VS$ given by
the form (\ref{227})
\[
\om\w\otimes d\si^\mu_c [\dr^c_\mu +\frac18\eta^{cb}\si^a_\mu
[\g_a,g_b]^B{}_Ay^A \dr_B].
\]


\begin{thebibliography}{ederf}

\bibitem{car} J. Cari\~nena, M. Crampin and L. Ibort, First
order multisymplectic formalism in filed theory, {\it Differential
Geometry and its Application}, {\bf 1}  (1991) 345.

\bibitem{got} M. Gotay, A multisymplectic framework for
classical field theory and the calculus of variations, in {\it
Mechanics, Analysis and Geometry: 200 Years after Lagrange}, ed.
M.Francaviglia (Elseiver Science Publishers B.V., 1991)
p. 203.

\bibitem{gun} C. G\"unther, The polisymplectic Hamiltonian
formalism in field theory and calculus of variations,  {\it Journal.
of Differential Geometry}, {\bf 25} (1987) 23.

\bibitem{kol} I.Kola\v{r}, P.W.Michor, J.Slov\'ak,
{\it Natural operations in differential geometry}, (Springer-Verlag, Berlin,
 1993).

\bibitem{2sar}  G. Sardanashvily, On the geometry of
spontaneous symmetry breaking, {\it
Journal of Mathematical Physics}, {\bf 33} (1992) 1546.

\bibitem{3sar} G. Sardanashvily  and O. Zakharov, {\it  Gauge
Gravitation Theory}  (World Scientific, Singapore, 1992).

\bibitem{6sar}
G. Sardanashvily and O. Zakharov,
On application of the Hamilton
formalism in fibred manifolds to field theory, {\it
Differential Geometry
and its Applications},  {\bf 3}, (1993), 245.

\bibitem{sard} G.Sardanashvily, {\it Gauge Theory in Jet Manifolds}
(Hadronic Press, Palm Harbor, 1993).

\bibitem{lsar}  G. Sardanashvily,
Constraint field systems in multimomentum canonical
variables, {\it  Journal of Mathematical Physics},
{\bf 35} (1994) 6584.

\bibitem{10sar}
G.Sardanashvily, Gravity as a Higgs Field, E-prints: gr-qc/9405013,
9407032, 9411013.

\bibitem{9sar}
G.Sardanashvily, Hamiltonian field systems on composite manifolds,
E-print: hep-th/9409159.

\bibitem{8sar} G.Sardanashvily, Five lectures on the jet manifold methods in
field theory, E-print: hep-th/9411089.

\bibitem{sau} D. Saunders, {\it The Geometry of Jet Bundles},
L.M.S. Lect. Notes Ser., {\bf 142}
(Cambridge Univ. Press, Cambridge, 1989).

\end{thebibliography}
\end{document}